\documentclass{article}
\usepackage{graphicx}
\usepackage{onecolceurws}
\usepackage{amsmath}
\usepackage[colorlinks=true,citecolor=blue,linkcolor=blue,urlcolor=blue]{hyperref}
\usepackage{nicefrac}

\usepackage[utf8x]{inputenc}
\usepackage{amssymb, upgreek}

\usepackage{color}
\definecolor{keywordcolor}{rgb}{0.7, 0.1, 0.1}   
\definecolor{commentcolor}{rgb}{0.4, 0.4, 0.4}   
\definecolor{symbolcolor}{rgb}{0.0, 0.1, 0.6}    
\definecolor{sortcolor}{rgb}{0.1, 0.5, 0.1}      
\definecolor{errorcolor}{rgb}{1, 0, 0}           
\definecolor{stringcolor}{rgb}{0.5, 0.3, 0.2}    

\usepackage{listings}

\newcommand{\bR}{\mathbb R}

\title{Verified Optimization\\(work in progress)}

\author{
Alexander Bentkamp \\ Vrije Universiteit Amsterdam\\
Department of Computer Science \\
De Boelelaan 1111, \\1081 HV Amsterdam,  Netherlands \\ bentkamp@gmail.com
\and
Jeremy Avigad \\
Carnegie Mellon University\\
Department of Philosophy\\
Baker Hall 161, \\Pittsburgh, PA 15213\\ avigad@cmu.edu
}

\institution{}

\begin{document}
\maketitle

\begin{abstract}
Optimization is used extensively in engineering, industry, and finance,
and various methods are used to transform problems to the point
where they are amenable to solution by numerical methods.
We describe progress towards developing a framework, based on the Lean
interactive proof assistant,
for designing and applying such reductions in reliable and flexible ways.
\end{abstract}

\section{Introduction}
\label{section:introduction}

Interactive proof assistants are used to verify complex mathematical claims
with respect to the primitives and rules of a formal axiomatic foundation.
Formalization yields a high degree of certainty in the correctness of such claims,
but it places a very high burden on practitioners, and
for many purposes it is a higher standard than users may want or need.
The project we describe here is motivated by the observation that
interactive theorem provers can offer a wider range of benefits
to applied mathematicians.
Sometimes even just a formal specification of a complex problem or model
is helpful, since it provides clarity and precision that can
serve as a touchstone for informal reasoning and algorithmic implementation.
Formal representation in a theorem prover can also serve as a gateway
to the use of external tools like computer algebra systems,
numeric computation packages,
and automated reasoning systems,
providing a basis for coordinating and interpreting the results.
And verification itself is not an all-or-nothing affair;
every working mathematician and scientist has to balance pragmatic constraints
against the goal of ensuring that their results are as reliable as possible,
and they should have the flexibility of deciding where verification effort
matters the most.
Proof assistants need to become a help rather than a hindrance before they
are ready to enter the mainstream.

Optimization problems and constraint satisfaction problems are now
ubiquitous in engineering, industry, and finance.
These address the problem of finding an element of $\bR^n$ satisfying
a finite set of constraints
or determining that the constraints are unsatisfiable;
the problem of bounding the value of an objective function over the domain
defined by such a set of constraints;
and the problem of finding a value of the domain that maximizes (or minimizes)
the value of the objective function.
Linear programming,
revolutionized by Dantzig's introduction of the simplex algorithm in 1947,
deals with the case in which the constraints and objective function are linear.
The development of interior point methods in the 1980s
allows for the efficient solution of problems defined by convex constraints
and objective functions, which gives rise to the field of convex programming \cite{boyd-vandenberghe-2014-convex}.

There are a number of ways in which formal verification can be used to
improve the reliability of such methods. Checking the correctness of a solution
to a satisfaction problem is easy in principle: one simply plugs the result into
the constraints and checks that they hold.
Verifying the correctness of a bounding problem or optimization problem is
often almost as easy, in principle, since the results are often underwritten
by the existence of suitable \emph{certificates}
that are output by the optimization tools.
In practice, these tasks are made more difficult
by the fact that floating point calculation can introduce
numerical errors that bear on the correctness of the solution.

Here we focus on a different part of the process,
namely, that of manipulating a problem and reducing it to a form
where optimization software can be applied.
Mathematical models are often complex, and practitioners rely on heuristics
and expertise to put problems into forms that admit computational solutions.
Such transformations are hard to automate,
and manual transformation is error-prone.
Our goal is to show proof assistants can be put to good use towards
finding and verifying these transformations,
and to develop tools to support the process.

In Section~\ref{section:optimization:problems:and:reductions},
we describe a general formal framework for reasoning about problems
and reductions between them.
In Section~\ref{section:disciplined:convex:programming},
to illustrate the method,
we describe a class of transformations
that form the basis for \emph{disciplined convex programming},
a component of the popular CVX package \cite{grant-boyd-2014-cvx}.
Even though these transformations are relatively straightforward,
we argue that a proof assistant provides a natural setting to
carry them out in a verified way.
In Section~\ref{section:future:plans}, we discuss more substantial problem transformations and reductions.

Our current work is spread between two versions of the Lean system \cite{de-moura-et-al-2015-lean}.
The development of Lean 3 has mostly stabilized; its library, {\sf mathlib} \cite{mathlib-2020},
comprises around 600\,000 lines of code 
and covers substantial portions of algebra, linear algebra, topology, measure theory, and analysis. Lean 4 is currently under development
as a performant dependently typed programming language;
it is not backward compatible with Lean 3, and does not have a substantial library yet.
We intend to make use of Lean 4's support for user extensible syntax,
as described below, and we plan to move the full development to Lean 4 as
soon as its library will support it.

\section{Optimization Problems and Reductions}
\label{section:optimization:problems:and:reductions}

The general structure of an optimization problem is as follows:
\begin{lstlisting}
  structure Minimization :=
  (Domain : Type)
  (objFun : Domain → ℝ)
  (constraints : Domain → Prop)
\end{lstlisting}
We express maximization problems by negating the objective function.
We assume that the objective function \lstinline{objFun} is defined over the data type
\lstinline{Domain} and takes values in the real numbers.
(It is often useful to allow values in the extended real numbers,
and we have not ruled out adopting this option instead.)
The domain is often $\bR^n$ or a space of matrices,
but it can also be something more exotic,
like a space of functions.
A \emph{feasible point} is an element of the domain satisfying the constraints.
A \emph{solution} to the minimization problem is a feasible point \lstinline{x} such that
for every feasible point \lstinline{y}
the value of the objective function at \lstinline{x} is smaller than
or equal to the value at \lstinline{y}.

Feasibility and bounding problems can also be expressed in these terms.
If the objective function is constant (for example, the constant zero function),
a solution to the optimization problem is simply a feasible point.
And given a domain, an objective function, and constraints,
a value \lstinline{b} is a (strict) bound on the value of the objective function over the domain
if and only if the feasibility problem obtained by adding the inequality \lstinline{objFun x ≤ b} to the constraints has no solution.

If \lstinline{p} and \lstinline{q} are problems, a \emph{reduction} from \lstinline{p} to \lstinline{q} is simply a function
mapping any solution to \lstinline{q} to a solution to \lstinline{p}.
The existence of such a reduction means that to solve \lstinline{p} it suffices to solve \lstinline{q}.
If \lstinline{p} is a feasibility problem, it means that the feasibility of \lstinline{q} implies the
feasibility of \lstinline{p}, and, conversely,
that the infeasibility of \lstinline{p} implies the infeasibility of \lstinline{q}.
With this framework in place, we can now easily describe what we are after:
we are looking for a system that helps a user reduce a problem \lstinline{p} to
a problem \lstinline{q} that can be solved by an external solver.
(For a bounding problem \lstinline{q}, the goal is instead to find a reduction to \lstinline{q} from an infeasible problem \lstinline{p}.)
At the same time, we wish to verify the correctness of the reduction,
either automatically or with user interaction.
This will ensure that the results from the external solver really address
the problem that the user is interested in solving.

There are at least three advantages to performing such transformations in a proof assistant.
First, it offers strong guarantees that the results are correct and have
the intended meaning.
Second, it means that users can perform the transformations interactively or partially,
and thus introspect and explore the results of individual transformation steps.
Finally, users can benefit from the ambient mathematical library,
including a database of functions and their properties.

The formulation described above is good for reasoning about problems in general,
but it is not as good for reasoning about \emph{particular} problems.
The optimization function in our representation of a minimization problem
is a unary function and the constraints are given by a unary predicate,
but we commonly think of these in terms of multiple variables, along these lines:
\begin{lstlisting}
minimization
  !vars x y z
  !objective x + y + z
  !constraints
    x + 3 * y > 5,
    z < 10
\end{lstlisting}
We have implemented exactly this syntax using Lean 4's flexible
mechanisms for macro expansion \cite{ullrich-de-mourra-2020-beyond-notations},
so that it represents an expression of the \lstinline{Minimization} type
presented above.
We use exclamation marks because the keywords we choose
become parser tokens in any Lean file that imports our library;
for example, with the syntax above, the tokens \lstinline{minimization},
\lstinline{!vars}, \lstinline{!objectives}, and \lstinline{!constraints}
can no longer be used as variable names or identifiers.
Conversely, Lean's \emph{delaborator} makes it possible
to pretty-print suitably-described problems in the form above.
We intend to use tactics written in Lean's powerful
metaprogramming language and a supporting library
to facilitate the interactive construction of reductions,
using the means above to mediate between internal representations
and user-facing syntax.


\section{Disciplined Convex Programming}
\label{section:disciplined:convex:programming}

Disciplined convex programming (DCP) \cite{grant-2004-thesis,grant-et-al-2006-dcp}
is a framework to specify convex optimization problems.
Any optimization problem following the rules of the framework can
be solved fully automatically.
An example of a DCP problem is the following \cite[equation (3.60)]{grant-2004-thesis}:
\begin{align*}
\text{minimize}\quad cx\tag{$\star$}\quad
\text{subject to}\quad&{\exp}(y)\le\log(a\sqrt x + b),
\quad ax + by = d
\end{align*}
where $a,b,c,d$ are parameters with $a \geq 0$, and $x,y$ are variables.
For this problem to be DCP conformant, it is crucial
that for instance the argument $a\sqrt x + b$ of the concave, nondecreasing function $\log$ is itself concave;
that the convex expression $\exp(y)$ is on the left of $\le$; and that
the concave expression $\log(a\sqrt x + b)$ is on the right of~$\le$.
The DCP rules allow for a systematic verification of all necessary conditions.

The DCP framework is implemented in the modeling systems
CVX \cite{grant-boyd-2014-cvx,grant-boyd-2008-cvx},
CVXPy \cite{agrawal-et-al-2018-cvxpy,diamond-boyd-2016-cvxpy},
Convex.jl \cite{udell-et-al-2014-convexjl}, and
CVXR \cite{fu-2020-cvxr}; a related reduction system is implemented in YALMIP~\cite{lofberg-2004-yalmip}.
These systems transform a DCP problem
into conic form, a more constrained canonical form that subsumes
linear, quadratic, and semidefinite programs.
The conic problem is then solved by external solvers such as SeDuMi \cite{sturm-1999-sedumi} and
SDPT3 \cite{toh-1999-sdpt3},
and the result is translated back into a solution to the original DCP problem.

Conic problems have the following form,
where $c \in \mathbb{R}^n$,
$A \in \mathbb{R}^{m\times n}$,
$b \in \mathbb{R}^m$,
$G \in \mathbb{R}^{k\times n}$,
$h \in \mathbb{R}^k$
are parameters
and $x \in \mathbb{R}^n$ are variables.
\begin{align*}
\text{minimize}\quad c^t x\quad
\text{subject to}\quad&Ax = b,\quad
Gx - h \in K
\end{align*}
The parameter $K$ is a convex cone. The conic solvers require $K$ to be a cartesian product of cones supported by the solver---e.g.,
the nonnegative orthant $\{x \in \mathbb{R}^l \mid x_i \geq 0 \text{ for all $i$}\}$,
the second-order cone $\{x \in \mathbb{R}^l \mid x_1 \geq \sqrt{x_2^2 + \cdots x_l^2}\}$,
or the exponential cone $\{x \in \mathbb{R}^3 \mid x_1 \geq x_2 e^{\nicefrac{x_3}{x_2}}, x_2 > 0\} \cup \{x \in \mathbb{R}^3 \mid x_1 \geq 0, x_2 = 0, x_3 \leq 0\}$.

Support for DCP in Lean will form the basis of our project.
DCP is widely applicable, and the transformations are relatively simple.
It is therefore a suitable testbed for the basic definitions described
in Section \ref{section:optimization:problems:and:reductions}
and for the potential of optimization tooling in proof assistants.

Our tool will accept a DCP problem specified in Lean and translate it into conic form while verifying that any solution of the conic
problem yields a solution of the original problem. In this paper, we will focus on this first task.
We envision that the tool will then send the conic problem to an external conic optimization solver.
The solver will return a solution, along with a dual solution that will allow us to verify the correctness of the
result independently in Lean.

We have considered sending the DCP problem directly to CVX or a similar high-level modeling system.
However, to the best of our knowledge, the dual solutions that CVX provides do not allow us
to verify the result independently, without canonizing the problem to conic form in Lean.
Moreover, experimenting with these transformations will help us to prepare the groundwork
for more complex problem transformations that cannot be handled fully automatically.

We demonstrate a problem transformation in Lean using the DCP program $(\star)$ above.
The transformation is shown in Figure~\ref{figure:transformation}.
Ultimately, we would like to fully automate such DCP transformations and to require user interaction
only for more complex reductions.
The program $(\star)$ is formulated above as found in Grant's thesis \cite{grant-2004-thesis}.
Grant assumes that $\log(x)$ for $x \le 0$ and $\sqrt{x}$ for $x < 0$ take value $- \infty$.
We could do the same in Lean, but algebraically the extended real numbers are not a convenient
number system to work with and Lean's library is substantially more comprehensive for the reals.
Instead, we add the constraints \lstinline{0 ≤ x} and \lstinline{0 < a * sqrt x + b} explicitly,
resulting in the definition of \lstinline{prob₁}.

\begin{figure}
\begin{minipage}[t]{.33\textwidth}
\begin{lstlisting}
def prob₁ := minimization
  !vars x y
  !objective c * x
  !constraints
    exp y
      ≤ log (a * sqrt x + b),
    a * x + b * y = d,
    0 ≤ x,
    0 < a * sqrt x + b
\end{lstlisting}
\end{minipage}%
\begin{minipage}[t]{.33\textwidth}
\begin{lstlisting}
def prob₂ := minimization
  !vars t₁ x y
  !objective c * x
  !constraints
    exp y ≤ t₁,
    t₁ ≤ log (a * sqrt x + b),
    a * x + b * y = d,
    0 ≤ x,
    0 < a * sqrt x + b
\end{lstlisting}
\end{minipage}%
\begin{minipage}[t]{.33\textwidth}
\begin{lstlisting}
def prob₃ := minimization
  !vars t₂ t₁ x y
  !objective c * x
  !constraints
    t₂ ^ 2 ≤ x,
    exp y ≤ t₁,
    t₁ ≤ log (a * t₂ + b),
    a * x + b * y = d,
    0 ≤ x,
    0 < a * t₂ + b
\end{lstlisting}
\end{minipage}%

\begin{minipage}[t]{.33\textwidth}
\begin{lstlisting}
def prob₄ := minimization
  !vars t₃ t₂ t₁ x y
  !objective c * x
  !constraints
    exp t₃ ≤ a * t₂ + b,
    t₂ ^ 2 ≤ x,
    exp y ≤ t₁,
    t₁ ≤ t₃,
    a * x + b * y = d,
    0 ≤ x,
    0 < a * t₂ + b
\end{lstlisting}
\end{minipage}%
\begin{minipage}[t]{.33\textwidth}
\begin{lstlisting}
def prob₅ := minimization
  !vars t₃ t₂ t₁ x y
  !objective c * x
  !constraints
    exp t₃ ≤ a * t₂ + b,
    t₂ ^ 2 ≤ x,
    exp y ≤ t₁,
    t₁ ≤ t₃,
    a * x + b * y = d
\end{lstlisting}
\end{minipage}%
\begin{minipage}[t]{.33\textwidth}
\end{minipage}%
\vspace{-\baselineskip}
\caption{A DCP transformation in Lean}
\label{figure:transformation}
\end{figure}

\begin{figure}
\begin{minipage}[t]{.5\textwidth}
\begin{lstlisting}
def linearization_antimono (D : Type)
(f g : D → ℝ) (c : ℝ → D → Prop)
(h_mono: ∀ x r s, r ≤ s → c s x → c r x) :
Minimization.Reduction
{ Domain := D,
  objFun := f,
  constraints := fun x => c (g x) x }
{ Domain := ℝ × D,
  objFun := fun y => f y.2,
  constraints :=
    fun y => g y.2 ≤ y.1 ∧ c y.1 y.2 }
\end{lstlisting}
\end{minipage}%
\begin{minipage}[t]{.5\textwidth}
\begin{lstlisting}
def graph_expansion_concave (D : Type)
(f g : D → ℝ) (c d : ℝ → D → Prop)
(h_mono: ∀ x r s, r ≤ s → c r x → c s x)
(hg : ∀ x v, c v x
    → isGreatest {y | d y x} (g x)) :
Minimization.Reduction
{ Domain := D,
  objFun := f,
  constraints := fun x => c (g x) x }
{ Domain := ℝ × D,
  objFun := fun y => f y.2,
  constraints :=
    fun y => d y.1 y.2 ∧ c y.1 y.2 }
\end{lstlisting}
\end{minipage}
\caption{Examples of reduction schemas}
\label{figure:lemmas}
\end{figure}

The first transformation step, from \lstinline{prob₁} to \lstinline{prob₂},
consists of moving \lstinline{exp y} into a separate constraint.
Grant calls this kind of transformation \emph{linearization}.
The occurrence of \lstinline{exp y} is replaced by an auxiliary variable \lstinline{t₁}
and we add the constraint \lstinline{exp y ≤ t₁}.
The transformation can be justified by the schema \lstinline{linearization_antimono} in Figure~\ref{figure:lemmas},
which is parameterized by a function \lstinline{c} that is monotone in
its first argument.
Instantiating it appropriately yields a reduction from \lstinline{prob₁} to \lstinline{prob₂}.
The condition for the transformation is that \lstinline{exp y} occurs in an antimonotone
context---i.e., if the constraints hold for some value \lstinline{s} of \lstinline{exp y},
then they also hold for values smaller than \lstinline{s}.
There is an analogous reduction schema \lstinline{linearization_mono} for monotone contexts that
would introduce the constraint \lstinline{t₁ ≤ exp y} instead.

Second, to reduce \lstinline{prob₂} to \lstinline{prob₃}, we eliminate the occurrence of \lstinline{sqrt}
by replacing it by its \emph{graph implementation}.
A graph implementation is a description of a concave [convex] function
as a convex maximization [minimization] problem.
For example, for $x \geq 0$, $\sqrt x$ can be described as the greatest number $y$ such that  $y ^ 2 ≤ x$.
The reduction schema \lstinline{graph_expansion_concave} in Figure~\ref{figure:lemmas} justifies the process of replacing a concave function
by its graph implementation, called \emph{graph expansion}.
In our case, we replace both occurrences of \lstinline{sqrt x} by an auxiliary variable \lstinline{t₂}
and add the constraint \lstinline{t₂ ^ 2 ≤ x},
which yields a reduction from \lstinline{prob₂} to \lstinline{prob₃}.
As for linearization, the condition for graph expansion is that \lstinline{sqrt x}
occurs in a monotone context. In fact, linearization is a special case
of graph expansion using the trivial graph implementation
defining $f(x)$ as the greatest number $y$ such that $y \leq f(x)$.
Next, we perform a graph expansion on \lstinline{log (a * t₂ + b)}, yielding \lstinline{prob₄} and a reduction from \lstinline{prob₃} to \lstinline{prob₄}.

Once \lstinline{log} and \lstinline{sqrt} have been eliminated, the constraints \lstinline{0 ≤ x} and
\lstinline{0 < a * t₂ + b} have served their purpose and can be removed, yielding \lstinline{prob₅}.
Using the constraints \lstinline{exp t₃ ≤ a * t₂ + b} and
\lstinline{t₂ ^ 2 ≤ x}, it is easy to show that these constraints are redundant.
For this step, we can even show that \lstinline{prob₄ = prob₅}.

Finally, problem \lstinline{prob₅} can be written in conic form.
The constraint \lstinline{a * x + b * y = d} constitutes the linear component of the conic form.
We can write
\lstinline{t₁ ≤ t₃} as a nonnegative orthant constraint;
\lstinline{t₂ ^ 2 ≤ x} as a second-order cone constraint;
\lstinline{exp t₃ ≤ a * t₂ + b} and
\lstinline{exp y ≤ t₁} as an exponential cone constraint.
The product of these cones constitutes the cone $K$ of the conic form.

We have proven the above transformation correct in Lean by
applying the reductions manually and proving the side conditions.
Our goal is to fully automate DCP canonization.
To this end, we will need a
library of graph implementations,
a tactic for proving monotonicity,
a tactic to derive the actual conic form from the fully graph expanded problem,
and an overarching tactic to guide the transformation process.


\section{Future Plans}
\label{section:future:plans}

Our approach will reveal its full potential when dealing
with problem transformations that are hard to automate because they are
specific to a particular problem.
The DCP methodology relies on experts to program the necessary graph implementations into the system.
We believe that a verified toolbox makes it easier for users---both experts and novices---to extend the system
to handle new reductions and to get the details right.
For extensions of the DCP methodology such as
disciplined convex-concave programming \cite{shen-et-al-2016-dccp},
disciplined geometric programming \cite{agrawal-et-al-2019-dgp}, and
disciplined quasiconvex programming \cite{agrawal-boyd-2020-dqp},
this is even more crucial.
For instance,
disciplined convex-concave programming requires information about sub- and supergradients, and
quasiconvex programming requires representations of the sublevel sets of the employed quasi-convex functions.
Another example is the barrier method that requires the user to come up with appropriate barrier functions.
It is impossible to devise a library that includes all functions that users will ever need,
but the verified approach provides a safe environment to derive the required information interactively.



We aim to test and evaluate our toolset with concrete applications.
Optimization and feasibility tools are often used in control theory to establish
stability and asymptotic stability of systems, as well as safety properties.
There is now a substantial literature on the use of formal methods to support this, and,
in particular, to develop ways of reliably reducing verification problems to problems
that can be checked by symbolic and numeric methods.
We believe that a system like the one we are developing can contribute in two ways:
first, by providing a general mathematical library and tools to verify the soundness of the theoretical reductions,
and second, by providing an interactive tool for \emph{applying} the reductions to specific problems,
ensuring that the data is in the right form and that the side conditions are met.

For example, a recent paper by Wang et al.~\cite{wang-et-al-2021-invariant-barrier} that we might explore in a case study
uses optimization to synthesize invariants of hybrid systems and thereby prove safety over the infinite time horizon.
The synthesized invariant is a barrier certificate,
encoded as an optimization problem constrained by bilinear matrix inequalities.
To make the problem amenable to conic solvers, Wang et al.\ transform the inequalities
into difference-of-convex constraints.
Through a transformation resembling disciplined convex-concave programming,
they bring the problem into conic form.
Finally, they use the branch-and-bound framework to ensure finding the global optimum, not a local one.
We consider these transformations to be an excellent case study for our approach because they
are practical and hard to automate.

\subsubsection*{Related work}

Formal methods have been used to solve bounding problems \cite{ratschan-07-safety-verification, gao-et-al-12-delta-complete}, constraint satisfaction problems \cite{franzle-et-al-07-constraint-systems}, and optimization problems \cite{kong-et-al-07-exists-forall}. The literature is too large to cover here; \cite{deshmukh-et-al-19-verification-cyber-physical} surveys some of the methods that are used in connection with the verification of cyber-physical systems.

Proof assistants have been used to verify bounds in various ways. Some approaches use certificates from numerical packages; Harrison \cite{harrison-07-sum-of-squares} uses certificates from semidefinite programming in HOL Light, and Magron et al.~\cite{magron:et:al:15} and Martin-Dorel and Roux \cite{dorel-roux-17-valid-sdp} use similar certificates in Coq.
Solovyev and Hales use a combination of symbolic and numeric methods in HOL Light~\cite{solovyev-hales-14-nonlinear-inequalities}. Other approaches have focused on verifying symbolic and numeric algorithms instead. For example, Mu{\~{n}}oz, Narkawicz, and
Dutle \cite{munoz-narkawicz-dutle-18-univariate-pvs} verify a decision procedure for univariate real arithmetic in PVS and Cordwell, Tan, and Platzer \cite{cordwell-tan-platzer-21-univariate} verify another one in Isabelle.
Narkawicz and Mu\~noz \cite{narkawicz-munoz-13-bnb} have devised a verified numeric algorithm to find bounds and global optima. Cohen et al.~\cite{cohen-davy-feron-garoche-17,cohen-feron-garoche-2020-verification-convex} have developed a framework for verifying optimization algorithms using the ANSI/ISO C Specification Language (ACSL) \cite{acsl-2020}.

Although the notion of a convex set has been formalized in a number of theorem provers, we do not know of any full development of convex analysis.
The Isabelle \cite{nipkow-paulson-wenzel-2002} {\sf HOL-Analysis} library\footnote{\url{https://isabelle.in.tum.de/dist/library/HOL/HOL-Analysis/}} includes properties
of convex sets and functions, including Carath\'eodory's theorem on convex hulls, Radon's theorem, and Helly's theorem, as well as properties of convex sets and functions on normed spaces and Euclidean spaces. A theory of lower semicontinuous functions by Grechuk \cite{grechuk-2011-lower-semicontinuous-afp} in the Archive of Formal Proofs \cite{blanchette-et-al-2015-afp} includes properties of convex functions. Lean's {\sf mathlib} includes a number of fundamental results,\footnote{\url{https://github.com/leanprover-community/mathlib/tree/master/src/analysis/convex}} including a formalization of the Riesz extension theorem by Kudryashov and Dupuis and a formalization of Jensen's inequality by Kudryashov. Allamigeon and Katz have formalized a theory of convex polyhedra in Coq with an eye towards applications to linear optimization \cite{allamigeon-katz-2019-formalization-convex-polyhedra}. We do not know of any project that has formalized reductions between optimization problems.


\subsubsection*{Acknowledgements} We are grateful to Seulkee Baek, Geir Dullerud, Paul Jackson, John Miller,  Ramon Fern\'andez Mir, Ivan Papusha, and Ufuk Topcu for
helpful discussions and advice. We also thank the anonymous reviewers for their corrections and suggestions.

\bibliographystyle{alpha}
\bibliography{bib}

\end{document}